\documentstyle[12pt]{article}
\parskip=\medskipamount                 % space between paragraphs
\thicklines                     % thick straight lines for pictures
\newcommand{\bimn}[7]{\bibitem{#1}#2,
{\em #3},
{ #4}$\;${\bf
#5}$\;$(#6)$\;${#7}.}
\def\inbar{\vrule height1.5ex width.4pt depth0pt}
\def\IC{\relax\,\hbox{$\inbar\kern-.3em{\rm C}$}}
\def\IN{\relax{\rm I\kern-.18em N}}
\def\IQ{\relax\,\hbox{$\inbar\kern-.3em{\rm Q}$}}
\def\IR{\relax{\rm I\kern-.18em R}}
\def\ZZ{\relax{\sf Z\kern-.4em Z}}
\def\a{\alpha} \def\b{\beta}

 \def\cK{{\cal K}} \def\cL{{\cal L}} 
 \def\cO{{\cal O}}  
 
\newtheorem{proposition}{Proposition}[section]

\newtheorem{conjecture}{Conjecture}[section]
\newtheorem{lemma}{Lemma}[section]
\newtheorem{definition}{Definition}[section]
\newtheorem{assumption}{Assumption}[section]

\catcode`\@=11

\newif\if@fewtab\@fewtabtrue

\catcode`\@=11

\newif\if@fewtab\@fewtabtrue

{\count255=\time\divide\count255 by 60
\xdef\hourmin{\number\count255}
\multiply\count255 by-60\advance\count255 by\time
\xdef\hourmin{\hourmin:\ifnum\count255<10 0\fi\the\count255}}
\def\ps@draft{\let\@mkboth\@gobbletwo
    \def\@oddhead{}
    \def\@oddfoot
       {\hbox to 7 cm{$\scriptstyle Draft\ version:\ \draftdate$
       \hfil}\hskip -7cm\hfil\rm\thepage \hfil}
    \def\@evenhead{}\let\@evenfoot\@oddfoot}

\def\ceqno{\global\@fewtabfalse
    \ifcase\@eqcnt \def\@tempa{& & &}\or \def\@tempa{& &}
      \or \def\@tempa{&}
      \or\def\@tempa{}\fi\@tempa
{\rm(\theequation)}}

\def\aeqno#1{\global\@fewtabfalse
    \ifcase\@eqcnt \def\@tempa{& & &}\or \def\@tempa{& &}
      \or \def\@tempa{&}
      \or\def\@tempa{}\fi\@tempa
{\rm(\theequation,#1)}}

\def\label#1{\ifnum\draftcontrol=1
 \global\def\draftnote{$\scriptstyle #1$}\fi
 \@bsphack\if@filesw {\let\thepage\relax
   \def\protect{\noexpand\noexpand\noexpand}%
\xdef\@gtempa{\write\@auxout{\string
      \newlabel{#1}{{\@currentlabel}{\thepage}}}}}\@gtempa
   \if@nobreak \ifvmode\nobreak\fi\fi\fi
  \@esphack}

\def\alabel#1#2{\label{#1}\global\@fewtabfalse
    \ifcase\@eqcnt \def\@tempa{& & &}\or \def\@tempa{& &}
      \or \def\@tempa{&}
      \or\def\@tempa{}\fi\@tempa
{\hbox to 3cm{\phantom{\rm(\theequation,#2)}
\draftnote \hfil}\hskip -3cm {\rm(\theequation,#2)}}}

\def\clabel#1{\label{#1}\global\@fewtabfalse
    \ifcase\@eqcnt \def\@tempa{& & &}\or \def\@tempa{& &}
      \or \def\@tempa{&}
      \or\def\@tempa{}\fi\@tempa
{\hbox to 3cm{\phantom{\rm(\theequation)}
\draftnote \hfil}\hskip -3cm{\rm(\theequation)}}}

\def\eqnarray{\def\draftnote{{}}\global\@fewtabtrue
\stepcounter{equation}\let\@currentlabel=\theequation
\global\@eqnswtrue
\global\@eqcnt\z@\tabskip\@centering\let\\=\@eqncr
$$\halign to \displaywidth\bgroup\@eqnsel\hskip\@centering\@eqcnt\z@
  $\displaystyle\tabskip\z@{##}$&\global\@eqcnt\@ne
   \hskip 1\arraycolsep \hfil$\displaystyle{##}$\hfil
  &\global\@eqcnt\tw@ \hskip 1\arraycolsep
$\displaystyle\tabskip\z@{##}$
\hfil  \tabskip\@centering&\global\@eqcnt\thr@@\llap{##}\tabskip\z@
\cr}

\def\endeqnarray{\@@eqncr\egroup
      \global\advance\c@equation\m@ne$$\global\@ignoretrue}

\def\@eqnnum{\hbox to 3cm{\phantom{\rm(\theequation)} \draftnote
                         \hfil}\hskip -3cm {\rm(\theequation)}}

\def\@@eqncr{\let\@tempa\relax
    \ifcase\@eqcnt \def\@tempa{& & &}\or \def\@tempa{& &}
      \or \def\@tempa{&}
      \or\def\@tempa{}
\fi\@tempa
\if@eqnsw
\if@fewtab\@eqnnum\fi
\stepcounter{equation}\fi\global
\@eqnswtrue\global\@eqcnt\z@\global\@fewtabtrue\cr}

\def\draftcite#1{\ifnum\draftcontrol=1#1\else{}\fi}

\def\@lbibitem[#1]#2{\item{}\hskip -3cm \hbox to 2cm
{\hfil$\scriptstyle\draftcite{#2}$}\hskip
1cm[\@biblabel{#1}]\if@filesw
     {\def\protect##1{\string ##1\space}\immediate
      \write\@auxout{\string\bibcite{#2}{#1}}}\fi\ignorespaces}

\def\@bibitem#1{\item\hskip -3cm \hbox to 2cm
{\hfil $\scriptstyle\draftcite{#1}$}\hskip 1cm
\if@filesw \immediate\write\@auxout
       {\string\bibcite{#1}{\the\value{\@listctr}}}\fi\ignorespaces}

\def\nsection#1{\section{#1}\setcounter{equation}{0}}

\def\nappendixe{\def\thesection{A}\section*{Appendix }
\def\theequation{{A.\arabic{equation}}}
\def\theproposition{{A.\arabic{proposition}}}
\setcounter{equation}{0}
\setcounter{proposition}{0}}

\def\draftdate{\number\month/\number\day/\number\year\ \ \ \hourmin }

\global\def\draftcontrol{0}
\catcode`\@=12

\def\theequation{{\thesection.\arabic{equation}}}

\def\qq{\begin{eqnarray}}
\def\qqq{\end{eqnarray}}
\def\rx#1{~(\ref{#1})}
\def\ex#1{eq.\rx{#1}}
\def\eex#1{eqs.\rx{#1}}
\def\cx#1{~\cite{#1}}
\def\rw#1{~\ref{#1}}

\newlength{\shiftwidth}
\def\shift#1{&&\hbox to \shiftwidth{\hfill $\displaystyle#1$}}
\newlength{\sshiftwidth}
\def\sshift#1{\lefteqn{\hbox to
\sshiftwidth{\hfill$\displaystyle#1$}}}
\def\snm{S_n(M)}

\def\ordH{| H_1(M,\ZZ)|}

\def\zmk{Z(M;k)}
\def\zsk{Z(S^3;k)}

\def\ztrmk{Z^{({\rm tr})}(M;k)}

\def\sjN{\sum_{j=1}^N}

\def\iptk{{i\pi\over 2K}}

\def\ipk{{i\pi\over K}}

\def\Dmn{D_{m,n}}

\def\spk{ \sin \left( {\pi\over K} \right) }

\def\ie{{\it i.e.\ }}
\def\eg{{\it e.g.\ }}
\def\cf{{\it cf.\ }}
\def\rhs{{\it r.h.s.\ }}
\def\lhs{{\it l.h.s.\ }}

\def\Rhs{RHS\ }
\def\wrt{WRT\ }

\def\asl{ASL\ }

\def\Tr{\mathop{{\rm Tr}}\nolimits}

\def\ord{\mathop{{\rm ord}}\nolimits}
\def\p{^{\prime}}

\def\prosign{\mathop{{\rm sign}}\nolimits}
\def\sign#1{\prosign\left(#1\right)}

\def\promod{\mathop{{\rm mod}}\nolimits}
\def\mod#1{\;(\promod #1)}

\def\pjn{ \prod_{j=1}^N }
\def\sjn{ \sum_{j=1}^N}

\def\va{ \vec{a} }
\def\vaa{ \vec{\a} }
\def\vv{ \vec{v} }

\def\ordH{\ord{|H_1(M,\ZZ)|} }

\def\lmu{ l^{(\mu)} }
\def\zmk{ Z(M;K) }
\def\zsk{ Z(S^3;K) }
\def\ztrmk{ Z^{(\rm tr)}(M;K) }

\def\snm{S_n(M)}
\def\dnm{\Delta_n(M)}

\def\aan{ \a_1,\ldots,\a_N }
\def\aant{ \a_2,\ldots,\a_N }
\def\aaqt{ (\a_2^2,\ldots,\a_N^2) }
\def\aaq{ (\a_1^2, \ldots, \a_N^2) }

\def\vaan{ \vaa_1, \ldots, \vaa_N }
\def\vaant{ \vaa_2, \ldots, \vaa_N }
\def\van{ \va_1, \ldots, \va_N }
\def\vant{ \va_2, \ldots, \va_N }

\def\lm{L_m}
\def\pml{P_{m,l}}
\def\lmaa{ L_m (\vaan) }
\def\lmaat{ L_m^{(2)} (\vaant) }

\def\pmlaa{ P_{m,l} (\vaan) }
\def\pmlaat{ P^{(2)}_{m,l} (\vaant) }

\def\lma{ L_m (\van) }
\def\lmat{ L_m^{(2)} (\vant) }

\def\pmla{ P_{m,l} (\van) }
\def\pmlat{ P^{(2)}_{m,l} (\vant) }

\def\ptmla{ \tilde{P}_{m,l} (\van) }
\def\ptml{ \tilde{P}_{m,l} }
\def\ptmlo{ \tilde{P}^{(1)}_{m,l} }
\def\ptmloaa{ \tilde{P}^{(1)}_{m,l} (\vaan) }

\def\Dmn{ D_{m,n} }
\def\Dmnt{ D^{(2)}_{m,n} }
\def\dmno{ d^{(1)}_{m,n} }

\def\jtr{ J^{({\rm tr})} }
\def\jan{ J_{\aan} }
\def\janlk{ \jan(\cL;K) }
\def\jtrat{ J^{({\rm tr})}_{\aant} }
\def\jtratb{ \jtrat(M_2,\cL\setminus\cL_1;K) }
\def\jtrathb{
J^{({\rm tr})}_{\a_3,\ldots,\a_N}(M_3,
\cL\setminus\cL_1\setminus\cL_2;K)}
\def\jtw{ J^{({\rm twist})}_{\aan} }
\def\jmer{ J^{({\rm meridian})}_{\aan;\b} }
\def\jf{ J_{\aan|\a} }

\def\sqtk{ \sqrt{2\over K} }
\def\spka{ \sin \left( {\pi\over K} \a \right) }
\def\spkaj{ \sin \left( {\pi\over K} \a_j \right) }
\def\spk{ \sin \left( {\pi \over K} \right) }
\def\sinb{ \sin \left( {\pi\over K} \b \right) }
\def\sinab{ \sin \left( {\pi\over K} \a\b \right) }

\def\iptk{ {i\pi \over 2K} }
\def\ipkt{ {i\pi K \over 2} }
\def\ipk{ \left( {i\pi\over K} \right) }
\def\ki{K^{-1}}
\def\pjna{ \left( \pjn \a_j \right)}
\def\pjnat{ \left( \prod_{n=2}^N \a_j \right)}

\def\smt{ \sum_{m\geq 2} }
\def\sml{ \sum_{m,l\geq 0 \atop m+l\neq 0} }

\def\sman{ \sum_{1\leq \a \leq \a_1 + \cdots + \a_N} }

\def\iaaj{ \int_{|\vaa_j| = \a_j} }
\def\iaj{ \int_{|\va_j| = {\a_j\over K} } }
\def\paaj{ \left( \pjn {d^2 \vaa_j \over 4\pi \a_j} \right) }
\def\paj{ \left( \pjn {K\over 4\pi} {d^2 \va_j \over|\va_j|}\right)}
\def\paajt{ \left( \prod_{j=2}^N {d^2 \vaa_j \over 4\pi \a_j} \right)}

\def\ffr{ \phi_{\rm fr} }
\def\ffro{ \phi_{\rm fr}^{(1)} }

\def\fact{ \frac{\sqrt{2}\pi} {K^{3\over 2}\ordH^{3\over 2} } }

\def\sinfv#1{\frac{\sin(\pi|\va_{#1}|)} {\pi |\va_{#1}| } }
\def\intva{ \int d^3\va_1\cdots d^3\va_N }

\def\ita#1{ \int\limits_{0\leq t_1 < \cdots < t_{#1} < 1}
  \!\!\!dt_1 \cdots dt_{#1} }
\def\itan{ \ita{n} }
\def\itanj{ \ita{n_j} }

\def\ddja#1{ {\delta\over \delta J(x_j(t_1))} \cdots
  {\delta \over \delta J(x_j(t_{#1}) ) } }
\def\ddjpa#1{ \left( \ddja{#1} \right) }
\def\ddjpan{ \left( {\delta\over \delta J(x_j(t_1))},
  \ldots,
  {\delta \over \delta J(x_j(t_n)) }\right) }

\def\ddjpanj{ \ddjpa{n_j} }

\def\ekiy#1{ \left. \exp \left[ #1 K^{1-m-n} \int_{S^3} dy_1\cdots
  dy_m\, \gmny\, J(y_1) \cdots J(y_m) \right] \right|_{J=0} }

\def\smtnz{ \sum_{m\geq 2 \atop n\geq 0} }
\def\smtnzs{ \sum_{ { m\geq 2 \atop n\geq 0 } \atop m+n\leq n_0 +1} }

\def\gmny{ G_{m,n}(y_1,\ldots,y_m) }

\def\ytof{ {1\over 4\pi} \epsilon_{\mu\nu\rho} {y_2^\rho - y_1^\rho
  \over |y_2 - y_1|^3} }
\def\dtz{ {1\over 2\pi i} {\delta(t_2 - t_1) \over z_2 - z_1} }

\def\zzbt{ (z_1,\bar{z}_1,t_1;z_2,\bar{z}_2,t_2) }

\def\cten#1{ C^{(#1)} }

\def\vdot{ v_1,\ldots,v_n }

\begin{document}
%\draft

\begin{titlepage}
%\centerline{\hfill                 UMTG-183-95}
\centerline{\hfill                 q-alg/9511025}
\vfill
\begin{center}

{\large \bf
On Finite Type Invariants of Links and Rational Homology Spheres
Derived from the Jones Polynomial and Witten-Reshetikhin-Turaev
Invariant. } \\

\bigskip
\centerline{L. Rozansky
%\footnote{Work supported by the National Science Foundation
%under Grant No. PHY-92 09978.}
}

\centerline{\em School of Mathematics, Institute for Advanced Study}
\centerline{\em Princeton, NJ 08540, U.S.A.}
\centerline{{\em E-mail address: rozansky@math.ias.edu}}

\vfill
{\bf Abstract}

\end{center}
\begin{quotation}
We define an infinite set of invariants of rational homology spheres
by presenting a link surgery formula which expresses them in terms of
the derivatives of the colored Jones polynomial of the link. We study
the properties of this formula and prove its invariance under the
Kirby moves.
\end{quotation}
\vfill
\end{titlepage}

\pagebreak

%+++++++++++++++++++++++++++++++
\nsection{Introduction}
\label{s1}
%+++++++++++++++++++++++++++++++

It has been established a while ago by D.~Bar-Natan\cx{BN1},
J.~Birman, X.~S.~Lin\cx{BL} and P.~Melvin, H.~Morton\cx{MeMo} that
the derivatives of the Jones polynomial with respect to the variable
$1/K$ ($q=e^{2\pi i \over K}$) at $K\rightarrow \infty$ are Vassiliev
(\ie finite type) invariants of knots and links. D.~Bar-Natan\cx{BN2}
and M.~Kontsevich\cx{Ko1} showed that these derivatives are related
to the Feynman diagram calculations of the Jones polynomial in the
framework of the quantum Chern-Simons theory proposed by
E.~Witten\cx{Wi1}. Guided by the principle that what is good for
knots is good for 3d manifolds, one might look for the finite type
invariants of manifolds among the Feynman diagram contributions to
the large $K$ limit of their Witten-Reshetikhin-Turaev (\wrt)
invariant.

A definition of finite type invariants of integer homology spheres
was given by  T.~Ohtsuki\cx{Oh3} (S.~Garoufalidis gave an
alternative definition in\cx{Ga1}). He demonstrated that Casson's
invariant was of finite type of order 3.
The definition of finite type invariants was later extended\cx{GO} to
rational homology spheres (\Rhs) (see also\cx{Ro4}). A particularly
promising place to look for these finite type invariants is the $1/K$
expansion of the trivial connection contribution to the \wrt
invariant. We conjectured a surgery formula for this
contribution\cx{Ro1},\cx{Ro3},\cx{Ro4}. This formula generates
explicit surgery formulas for the individual coefficients of the
$1/K$ expansion. We showed at the physical level of rigor that the
second coefficient is proportional to the Casson-Walker invariant.

In this paper we take a mathematically rigorous approach by using the
surgery formula of\cx{Ro1},\cx{Ro3} as a definition of the generating
function of an infinite sequence of invariants $\snm$ of a \Rhs. We
explore some of their properties and demonstrate their invariance
under the meridian Kirby move. The rigorous proofs
of the propositions stated in this paper will be presented
elsewhere.

In accordance with the physical considerations of\cx{Ro1} we expect
the invariants $\snm$ to coincide with $(n+1)$-loop corrections
studied by S.~Axelrod and I.~Singer\cx{AS1},\cx{AS2}. However the
results of this paper are independent of this identification or of any
physical arguments
used in our previous papers. The proofs are based upon
Kontsevish's integral representation\cx{Ko1} of the $1/K$ expansion
of the Jones polynomial and on the fusion properties of the local
cables established by N.~Reshetikhin and V.~Turaev\cx{RT}.

The propositions presented in this paper validate the claims of our
previous papers\cx{Ro4} and \cx{Ro5}. We sum up this claims in the
last section of this paper.

%++++++++++++++++++++++++++++++++++
\nsection{The $K^{-1}$ Expansion of the Colored Jones
Polynomial}
\label{s2}
%+++++++++++++++++++++++++++++++++++

Let $\cL$ be a framed $N$-component link in $S^3$.
We assign the $\a_j$-dimensional representations of $SU(2)$ to its
components. The colored Jones polynomial $\janlk$ of $\cL$ is
normalized in such a way that
%
%  and $\janlk$ be
%  its colored Jones polynomial normalized in such a way that
%
$$J_\a(\mbox{unknot};K) = {\spka \over \spk},\qquad
J_{\a_1,\a_2}(\cL_1\#\cL_2;K) = J_{\a_1}(\cL_1;K)
J_{\a_2}(\cL_2;K)$$
($\cL_1\#\cL_2$ denotes here a disconnected sum of $\cL_1$ and
$\cL_2$)
and changing the framing (\ie self-linking number) of a component
$\cL_j$ of $\cL$ by one unit leads to the multiplication of the whole
polynomial by a factor
$\exp \left( \iptk (\a_j^2 - 1) \right)$.
P.~Melvin and H.~Morton proved\cx{MeMo} that the
polynomial $\janlk$ can be expanded in powers of colors and
$\ki$:
\qq
\janlk = \pjna \sum_{n \geq 0 \atop 0\leq m\leq n} \Dmn\aaq K^{-n},
\label{1.1}
\qqq
here $\Dmn$ are homogeneous polynomials of order $m$.
whose coefficients are finite type invariants of $\cL$ of order $n$.
%
%D.~Bar-Natan\cx{BN1} and J.~Birman, X-S.~Lin\cx{BL} showed that the
%polynomials $\Dmn$ are finite type invariants of $\cL$ of order $n$.
%

The expansion\rx{1.1} can be presented in an alternative form:
\begin{proposition}
\label{p1.1}
For any link $\cL\in S^3$ there exists a set of $SU(2)$-invariant
homogeneous polynomials $\lmaa$, $\pmlaa$, $m,l\geq 0$ of order $m$
on the Lie algebra $su(2)$ (we denote the elements of $su(2)$ as
vectors $\vaa$) such that the colored Jones polynomial can be
presented as a multiple integral over the co-adjoint orbits
(\ie spheres in $\IR^3$ of radii $\a_j$ in the case of $SU(2)$)
corresponding to the representations assigned to the link components:
\qq
\janlk & = & \iaaj \paaj\; \exp \left( {i\pi\over 2} \smt \lmaa
K^{1-m}
\right.
\label{1.2}\\
&&\left. \qquad\qquad\qquad
 + \sml \pmlaa K^{-l-m} \right).
\nonumber
\qqq
This equation should be understood in the following way: for any
$n_0 > 0$
\qq
\iaaj \paaj \exp \left( {i\pi\over 2} \sum_{m=2}^{n_0+1} \lmaa
K^{1-m} + \!\!\!\!
\sum_{l,m\geq 0 \atop 0<l+m\leq n_0} \pmlaa K^{-l-m} \right)
&&
\nonumber\\
&&
\hspace*{-4.5in}
= \pjna \sum_{0\leq n\leq n_0 \atop 0\leq m\leq n} \Dmn\aaq K^{-n}
+ \cO(K^{-n_0 - 1} ).
\label{1.3}
\qqq
\end{proposition}

The idea to present the Jones polynomial as an integral over the
coadjoint orbits corresponding to the colors was first proposed by
N.~Reshetikhin\cx{Re}. We proved the formula\rx{1.2} at the physical
level of rigor in\cx{Ro2}. In fact, that proof becomes rigorous as
soon as one uses Kontsevich's integral\cx{Ko1} instead of a generic
Chern-Simons perturbation theory. We sketch this procedure in
Appendix. The detailed proof will be presented in\cx{Ro9}.

The next proposition was also proved in\cx{Ro2} (see the comments in
Appendix and\cx{Ro9}).
\begin{proposition}
\label{p01.1}
The polynomials $\lmaa$ satisfy the following properties:
\begin{enumerate}
\item
\qq
L_2 = \sum_{i,j=1}^N l_{ij} \vaa_i \cdot \vaa_j,
\label{1.4}
\qqq
here $l_{ij}$ are the linking numbers of $\cL$.
\item
If for the three numbers
$1 \leq j_1 \leq j_2 \leq j_3 \leq N$
all linking numbers $l_{j_{m_1} j_{m_2}}$, $1 \leq m_1 < m_2 \leq 3$
are zero then the sum of monomials of $L_3(\vaan)$ which depend only
on vectors $\vaa_{j_m}$, $1 \leq m \leq 3$ is
\qq
-4\pi \lmu_{j_1 j_2 j_3} \vaa_{j_1} \cdot (\vaa_{j_2} \times
\vaa_{j_3} ).
\label{1.5}
\qqq
here $\lmu_{ijk}$ are triple Milnor linking numbers.
\item
If for four numbers
$1 \leq j_1 \leq j_2 \leq j_3 \leq j_4 \leq N$
all the off-diagonal linking and triple Milnor linking numbers are
zero then the sum of monomials of $L_4(\vaan)$ containing only the
vectors $\vaa_{j_m}$, $1 \leq m \leq 4$ is
\qq
{\pi^2 \over 3} \sum_{ 1\leq m_1,m_2,m_3,m_4 \leq N}
\left (\lmu_{j_{m_1}j_{m_2}j_{m_3}j_{m_4}} -
\lmu_{j_{m_3}j_{m_1}j_{m_2}j_{m_4}} \right)
(\vaa_{j_{m_1}} \times \vaa_{j_{m_2}} ) \cdot
(\vaa_{j_{m_3}} \times \vaa_{j_{m_4}} ) ,
\label{1.7}
\qqq
here $\lmu_{ijkl}$ are quartic Milnor linking numbers.
\item
A polynomial $\lmaa$ is a linear combination of tree monomials coming
from tree graphs with trivalent vertices and $m$ legs as described in
\cx{Ro2}. A tree monomial is formed from a tree graph by placing the
Lie algebra structure constants ($\epsilon_{\mu\nu\rho}$ for $su(2)$)
at trivalent vertices, the Killing scalar product ($\delta_{\mu\nu}$
for $su(2)$) at edges and Lie algebra elements
$\vaa_j$, $1\leq j\leq N$ at legs.
\end{enumerate}
\end{proposition}
A finite type nature of Minlor's linking numbers was observed by
D.~Bar-Natan\cx{BN3} and X-S.Lin\cx{Li}.

It is instructive to see what happens to \ex{1.2} when $\cL$ has only
one component, \ie when $\cL$ is a knot $\cK$. First, since all the
polynomials $L_m,P_{m,l}$ are $SU(2)$-invariant, the integral over
the only variable $\vaa$ becomes trivial:
$\int_{|\vaa|=\a} {d^2\vaa \over 4\pi \a} = \a$.
Second, the property 4 of the Proposition~\ref{p01.1} implies that
all the polynomials $L_m(\vaa)$, $m\geq 3$ are equal to zero due to
the antisymmetric structure constants $\epsilon_{\mu\nu\rho}$ in the
vertices of the tree graphs. As a result, we end up with the
expansion
\qq
e^{-\iptk l_{11} \a^2 } J_\a(\cK;K) = \a
\left( 1 + \sml D_{m,l} \a^{2m} K^{-2m-l} \right)
\label{1.07}
\qqq
with some coefficients $D_{m,l}$. This expansion is equivalent to the
Melvin-Morton bound\cx{MeMo} proved by D.~Bar-Natan and
S.~Garoufalidis\cx{BG} (for a simple ``path integral'' proof
see\cx{Ro1}).

The formula\rx{1.2} can be put into a more suggestive form if we
introduce new variables
\qq
\va_j = {\vaa_j \over K}.
\label{1.8}
\qqq
In these variables
\qq
\janlk & = & \iaj \paj\; \exp \left( \ipkt \smt \lma
\right.
\label{1.9}\\
&&\left. \qquad\qquad\qquad
 + \sml \pmla K^{-l} \right).
\nonumber
\qqq
The integral decomposes into a product of the rapidly oscillating
exponential
$$\exp\left(\ipkt \smt L_m\right)$$
whose exponent is
proportional to $K$ and a slowly varying preexponential factor
$$\exp \left( \sml P_{m,l}K^{-l} \right)$$
whose expansion contains only
negative powers of $K$. As a result, the integral\rx{1.9} can be
calculated in the stationary phase approximation when
$K\rightarrow\infty$. In\cx{Ro2},\cx{Ro3} we discuss at physical
level of rigor a relation between the stationary phase configurations
of the vectors $\va_j$ and the homomorphisms of the link group into
$SU(2)$. We also establish (at the same rigor level) a
Melvin-Morton\cx{MeMo} type relation between the contribution of
configurations when all vectors $\va_j$ are (anti-)parallel and the
Alexander polynomial of $\cL$.

%+++++++++++++++++++++++++++++++
\nsection{Perturbative Invariants}
\label{s3}
%+++++++++++++++++++++++++++++++

We will use \ex{1.9} in order to define ``perturbative'' invariants
of \Rhs. We recall that if a 3d manifold $M$ is constructed by a
surgery on a framed link $\cL\in S^3$ (we denote this as
$M = \chi_\cL(S^3)$) then its \wrt invariant $\zmk$ is given by the
surgery formula
\qq
\zmk = \zsk e^{i\ffr} \sum_{1\leq \aan \leq K-1}
\janlk \pjn \sqtk \spkaj.
\label{1.10}
\qqq
Here $\zsk$ is the \wrt invariant of $S^3$:
\qq
\zsk = \sqtk \spk
\label{1.010}
\qqq
and $\ffr$ is a ``framing correction'' which depends on the surgery
data:
\qq
\ffr = - {3\over 4}\pi {K-2 \over K} \sign{l_{ij}},
\label{1.11}
\qqq
$\sign{l_{ij}}$ is the signature of the linking matrix $l_{ij}$.

We are going to use the formula\rx{1.9} for the Jones polynomial
in \ex{1.10} and substitute a sum $\sum_{1\leq\aan\leq K-1}$ in that
equation by an integral $\int_{0}^{\infty} d\a_1\cdots d\a_N$. The
integrals over $d\a_j$ and $d^2 \va_j$ combine into 3d integrals
$\int d^3\va_j$ which we will calculate in a stationary phase
approximation.

\begin{definition}
\label{d1.1}
Let $M$ be a \Rhs constructed by a surgery on a framed link
$\cL\in S^3$.
We define an infinite sequence of its invariants $\dnm$ (or,
equivalently, $\snm$) by expressing their generating function
\qq
\ztrmk & = &
\fact \; \sum_{n=0}^{\infty} \dnm K^{-n}
\label{1.12}\\
& = &
\fact \;\exp \left[ \sum_{n=1}^{\infty}\snm \ipk^n \right]
\nonumber
\qqq
by the surgery formula
\qq
\displaystyle
\ztrmk & = &\zsk e^{i\ffr} \int_{[\va_j=0]}
\left(\pjn \left({K\over 2}\right)^{3\over 2} d^3 \va_j \;\sinfv{j}
\right)
\label{1.13}\\
\displaystyle
&&\qquad
\times \exp \left[ \ipkt \smt \lma + \sml \pmla K^{-l} \right].
\nonumber
\qqq
Here the symbol $[\va_j=0]$ means that we are taking only the
contribution of the stationary phase point $\va_j=0$ to the integral
in the stationary phase approximation.

\end{definition}
Equation\rx{1.13} should be understood in the following way. Consider
the formal power series expansion
\qq
\left( \pjn \sinfv{j} \right)\; \exp \left( \ipkt \sum_{m=3}^{\infty}
\lma + \sml \pmla K^{-l} \right)\hspace{0.8in}
&&
\label{1.14}
\\
&&\hspace*{-2in}
= \sum_{m\geq 0 \atop l\geq -{m\over 3} } \ptmla K^{-l},
\nonumber
\qqq
here $\ptmla$ are homogeneous $SU(2)$-invariant polynomials of degree
$m$. Then for any $n_0\geq 0$
\qq
\lefteqn{
\sum_{n=0}^{n_0} \dnm K^{-n}
= \exp \left[ \sum_{n=1}^{n_0} \snm \ipk^n \right] + \cO(K^{-n_0-1})
}&&
\label{1.15}
\\
&=&
{K\over \pi}\spk e^{i\ffr} |\det (l_{ij})|^{3\over 2}
\left({K\over 2}\right)^{ {3\over 2}N}
\nonumber\\
&&\times
\intva\; \exp \left( \ipkt
\sum_{1\leq i,j \leq N} \!\! l_{ij} \va_i \cdot \va_j \right)
\hspace*{-0.1in}\!
\sum_{0\leq m\leq 6n_0 \atop -{m\over3}\leq l\leq n_0 - {m\over 2} }
\hspace*{-0.1in}
\ptmla K^{-l} + \cO(K^{-n_0 - 1})
\nonumber
\qqq
(note that
$\ordH = |\det(l_{ij})|$ ). Since only a finite number of polynomials
$\ptml$ participate in the preexponential sum of the \rhs of \ex{1.15}
for a given $n_0$, the integral there is well defined.

%++++++++++++++++++++++++++++++++++++++++++
\nsection{A Step-by-Step Procedure}
\label{s4}
%+++++++++++++++++++++++++++++++++++++++++++
For a given link $\cL$ \ex{1.2} does not determine the polynomials
$\lm$ and $\pml$. In other words, different sets of polynomials
$\lm,\pml$ may lead to the same Jones polynomial through \ex{1.2}.
However the invariants $\dnm$ and $\snm$ depend actually only on the
derivatives of the coefficients in the $1/K$ expansion\rx{1.1} of the
Jones polynomial $\janlk$.
%++++++++++++++++++++
\begin{proposition}
\label{p1.2}
For a given $n_0 > 0$ the invariants $\Delta_{n_0}(M)$ and
$S_{n_0}(M)$ can be expressed in terms of the linking numbers
$l_{ij}$ of $\cL$ and the polynomials $\Dmn\aaq$ of the
expansion\rx{1.1} for which
\qq
n\leq 2^N n_0,
\label{1.015}
\qqq
here $N$ is the number of components of $\cL$.
\end{proposition}
To prove this proposition we will use an alternative method of
calculating $\ztrmk$ by integrating over the colors of $\cL$
step-by-step\footnote{I am thankful to D.~Thurston for suggesting to
try this approach.}. It follows from the property 4 of the
Proposition~\ref{p01.1} that the maximum power of a given vector
$\va_j$ in any monomial of the polynomial $\lmaa$ is $m-2$ (at least
two legs in a tree with at least one vertex should be assigned
to different link components in order to get a non-zero monomial).
As a result, an expansion of the colored Jones polynomial of a link
satisfies the Melvin-Morton bound in individual colors:
%++++++++++++++++++++++++
\begin{proposition}
\label{p01.2}
The $1/K$ expansion\rx{1.1} can be rewritten as
\qq
e^{ -\iptk l_{11} \a_1^2} \janlk =
\pjna \sum_{m,n \geq 0 \atop m\leq {n\over 2} }
\dmno\aaqt\, \a_1^{2m} K^{-n},
\label{01.1}
\qqq
here $\dmno$ are polynomials of a degree not greater than $n$.
\end{proposition}

We can use this fact in order to substitute an integral instead of a
sum in the Reshetikhin-Turaev surgery formula for the surgery on
$\cL_1$. Then by switching a sum over $\a_1$ to an integral we arrive
at the definition which is similar to the Definition~\ref{d1.1}:
%+++++++++++++++++++++++++++
\begin{definition}
\label{d1.2}
If $l_{11}\neq 0$, then we define an infinite
sequence of invariants of
a link $\cL\setminus \cL_1$ in a \Rhs $M_2 = \chi_{\cL_1}(S^3)$ by
expressing their generating function
\qq
\jtratb = \pjnat \sum_{m,n\geq 0 \atop m\leq n} \Dmnt\aaq K^{-n}
\label{01.2}
\qqq
(here $\Dmnt$ are homogeneous polynomials of degree $m$) by a
stationary phase contribution of the point $a_1=0$ to the integral
\qq
\displaystyle
\lefteqn{  \displaystyle
\jtratb = e^{i\ffro} \sqrt{2K} \int\limits_{0\atop [a_1=0]}^{\infty}
da_1 e^{\ipkt l_{11} a_1^2} \sin(\pi a_1)
}\hspace*{3in}
\nonumber\\
\displaystyle
& &
\times \left( e^{-\ipkt l_{11} a_1^2} J_{Ka_1,\aant} (\cL;K)
\right),
\label{01.3}
\qqq
here
\qq
\ffro = -{3\over 4}\pi {K-2\over K} \sign{l_{11}}.
\label{01.4}
\qqq
\end{definition}
%+++++++++++++++++++
Equation\rx{01.3} should be understood in the following way: for any
$n_0 > 0$
\qq
\lefteqn{
e^{i\ffro} \sqrt{2K} \int\limits_{0}^{\infty} da_1\, e^{\ipkt l_{11}
a_1^2} \sin(\pi a_1) \sum_{0\leq m\leq n_0\atop 2m\leq n\leq m+n_0}
\dmno\aaqt a_1^{2m} K^{2m-n}
}\hspace*{3in}
\label{01.5}
\\
&=& \pjnat \sum_{0\leq n\leq n_0 \atop 0\leq m\leq n} \Dmnt\aaqt
K^{-n}.
\nonumber
\qqq
\begin{conjecture}
${\zsk\over\ztrmk} \jtratb$ is the trivial connection contribution to
the colored Jones polynomial of the link $\cL\setminus \cL_1$ in the
\Rhs $M_2 = \chi_{\cL_1}(S^3)$.
\end{conjecture}

We want to relate $\jtratb$ to the integral in \ex{1.13}.
\begin{proposition}
\label{p1.3}
The same generating function\rx{01.2} can be expressed as an integral
over $d^3\va_1$:
\qq
\lefteqn{
\jtratb=\iaaj \left( \prod_{j=2}^N {K\over 4\pi} {d^2\va_j\over
|\va_j|} \right) \, e^{i\ffro} \left({K\over 2}\right)^{3\over 2}
}
\label{01.6}\\
&&
\times
\int\limits_{\left[ \va_1 = -{1\over l_{11}} \sum_{j=2}^N l_{1j} \va_j
\right]}
\hspace*{-0.5in}
d^3 \va_1 \sinfv{1}
\exp \left( \ipkt \smt \lma + \sml \pmla K^{-l} \right)
\nonumber
\qqq
\end{proposition}

There are two distinct ways of understanding \ex{01.6} because there
are two ways of calculating the stationary phase integrals in its
\rhs. The first way is to substitute the condition
$|\va_j| = {\a_j\over K} $, $2\leq j\leq N$ inside the integral over
$\va_1$. Consider the expansion
\qq
\lefteqn{
e^{-\iptk l_{11} \a_1^2}
\exp\left( {i\pi\over 2} \smt \lmaa K^{1-m} + \sml \pmlaa K^{-m-l}
\right)
}\hspace*{4in}
\label{01.7}
\\
&=& \sum_{m,l\geq 0} \ptmloaa K^{-m-l},
\nonumber
\qqq
here $\ptmlo$ are invariant polynomials of total order at most
$2(m+l)$ and of the homogeneous order $m$ in $\vaa_1$. The fact that
$l\geq 0$ in the \rhs of \ex{01.7} follows from the property 4 of the
Proposition~\ref{p01.1} which limits the maximum power of $\vaa_1$ in
the polynomials $\lmaa$. In view of expansion\rx{01.7}, \ex{01.6}
reads in the first approximation for any $n_0>0$:
\qq
\lefteqn{
\left(\prod_{j=2}^N \a_j\right) \sum_{0\leq n\leq n_0\atop 0\leq
m\leq n} \Dmnt \aaqt K^{-n} =
e^{i\ffro} \left( {K\over 2}\right)^{3\over 2}
\int d^3\va_1 \sinfv{1} e^{\ipkt l_{11} \va_1^2}
}\hspace*{1in}
\label{01.8}\\
&&\times
\iaaj \left( \prod_{j=2}^N {d^2\vaa_j\over 4\pi \a_j} \right)
\sum_{0\leq m\leq 2n_0\atop 0\leq l\leq n_0 - {m\over 2} }
\ptmlo(\va_1,\vaant) K^{-l} + \cO(K^{-n_0-1}).
\nonumber
\qqq
This equation is easy to prove. We split the integral over $\va_1$:
$\int d^3\va_1 = \int_0^\infty da_1\int_{|\va_1|=a_1}d^2\va_1$ and
then integrate over $d^2\va_1$ and $d^2\va_j$, $2\leq j\leq N$.
According to the Proposition~\ref{p1.1} and \ex{01.1} we recover all
polynomials $\dmno\aaqt \a_1^{2m}K^{-n}$, $0\leq m\leq n_0$,
$2m\leq n\leq m+n_0$ from the polynomials
$\ptmlo(\va_1,\vaant) K^{-l}$ if we set $a_1 = {\a_1\over K}$.
Therefore the subsequent integral over $da_1$ gives the same results
in \eex{01.5} and \rx{01.8}.

The second way of calculating the integral over $d^3\va_1$ in the
\rhs of \ex{01.6} is to pretend that $|\va_j|$, $2\leq j\leq N$ are
of order $1$ and impose the condition
$|\va_j|={\a_j\over K}$ only after the integral over $d^3\va_1$ is
already calculated. The following lemma results from a Feynman
diagram  calculation\footnote{This calculation has nothing to do with
path integrals and is absolutely rigorous.} of the 3d integral over
$\va_1$ in \ex{01.6}.
%++++++++++++++++++++++++++++
\begin{lemma}
\label{l1.2}
$\!\!\!\!\!$
There exists a set of $SU(2)$-invariant homogeneous polynomials
$\lmat$, $\pmlat$, $m,l\geq 0$ of order m such that
\qq
\lefteqn{
e^{i\ffro} \left( {K\over 2} \right)^{3\over 2}
\hspace*{-0.4in}
\int\limits_{[\va_1 = -{1\over l_{11}}
\sum\limits_{j=2}^N l_{1j}\va_j] }
\hspace*{-0.4in}
d^3\va_1\; \sinfv{1} \exp \left[ \ipkt \smt \lma + \sml \pmla K^{-l}
\right]
}
\hspace*{1.1in}
\nonumber\\
&=& \exp \left[ \ipkt \sml \lmat + \sum_{m,l\geq 0}
\pmlat K^{-l} \right].
\label{01.08}
\qqq
The polynomials $L_m^{(2)}$ share the property 4 of the polynomials
$L_m$ in the Proposition~\ref{01.1} and
\qq
L_2^{(2)}(\vant) = \sum_{2\leq i,j \leq N} l\p_{ij} \va_i\cdot\va_j,
\qquad l\p_{ij}= l_{ij} - {l_{1i} l_{1j}\over l_{11} }.
\label{01.11}
\qqq
\end{lemma}
%+++++++++++++++++++++++
Combining \eex{01.6} and \rx{01.08} we find that a corollary of this
lemma is the following
%+++++++++++++++++++++++++
\begin{proposition}
\label{p1.4}
The invariant $\jtratb$ can be presented as an integral over the
coadjoint orbits in the way similar to the Proposition~\ref{p1.1}:
\qq
\lefteqn{
\jtratb = \iaaj \paajt \exp \left[ {i\pi\over 2} \sml \lmaat K^{1-m}
\right.
}\hspace*{3in}
\nonumber\\
&&
\left.
+
\sum_{m,l\geq 0} \pmlaat K^{-l-m} \right].
\label{01.10}
\qqq
\end{proposition}
The interpretation of this equation is similar to \ex{1.3}.

As we see, the function $\jtratb$ shares all the main properties of
the original Jones polynomial $\janlk$. Therefore the expansion of
$e^{-\iptk l\p_{22}\a_2^2} \jtratb$ in powers of $1/K$ satisfies the
same Melvin-Morton bound on individual powers of $\a_2$ as the
expansion of $e^{-\iptk l_{11} \a_1^2} \janlk$ in the
Proposition~\ref{p01.2}. This allows us to define the next invariant
$\jtrathb$ for $M_3 = \chi_{\cL_2}(M_2)$
in the way similar to \ex{01.3}:
\qq
\lefteqn{
\jtrathb = e^{i\ffr^{(2)}} \sqrt{2K}
\int\limits_{0 \atop [a_2=0]}^\infty da_2\,
e^{\ipkt l\p_{22} a_2^2 } \sin(\pi a_2)
}\hspace*{3in}
\label{01.012}\\
&&
\times
\left( e^{-\ipkt l\p_{22} a_2^2}
J^{({\rm tr})}_{Ka_2,\a_3,\ldots,\a_N}(M_2,\cL\setminus\cL_1;K)
\right),
\nonumber
\qqq
\qq
\ffr^{(2)} = -{3\over 4}\pi {K-2\over K} \sign{l\p_{22}}.
\label{01.1012}
\qqq
Repeating this procedure of step-by-step 1d stationary phase
integrations over $da_j$ $N$ times we end up with
$\jtr(M_{N+1};K)$ such that
\qq
\ztrmk = \zsk \jtr(M_{N+1};K).
%J^{({\rm tr})} (M;K) = {\ztrmk \over \zsk}.
\label{01.2012}
\qqq

It follows from \ex{01.5} that if we want to determine all the
polynomials $\Dmnt$ of \ex{01.2} for $n\leq n_0$ then it is enough to
know all the polynomials $\dmno$ of \ex{01.1} for $n\leq 2n_0$. This
means that at each step of consecutive stationary phase integrations
we reduce the precision in terms of powers of $K^{-1}$ by a factor of
2. Therefore in order
to calculate all the invariants $S_n(\chi_\cL(S^3))$,
$n\leq n_0$ by the step-by-step procedure it is enough to know all
the polynomials $\Dmn$ of \ex{1.1} for $n\leq 2^N n_0$, $N$ being the
number of components in $\cL$. This estimate may be an overkill, but
it almost proves the essential claim of the Proposition~\ref{p1.2}
that one has to know only a finite number of coefficients in the
expansion of the Jones polynomial in order to find $\ztrmk$ with any
given precision.

To complete the proof of the Proposition~\ref{p1.2} we have to decide
what to do if one of the self-linking numbers that appear in the
exponents of stationary phase integration formulas ($l_{11}$ in
\ex{01.3} or $l_{22}\p$ in \ex{01.012}) is zero so that the
stationary phase approximation calculation fails. We suggest to use a
``regularized'' linking matrix depending on parameters $\epsilon_j$
\qq
l^{(\epsilon)}_{ij} = l_{ij} + \epsilon_j \delta_{ij}
\label{01.14}
\qqq
instead of the actual linking matrix $l_{ij}$ of $\cL$ in all the
stationary phase calculations. Since we changed only the diagonal
part of $l_{ij}$, the ``regularized'' Jones polynomial is simply
related to the original one:
\qq
J^{(\epsilon)}_{\aan}(\cL;K) = e^{\iptk \sjn \epsilon_j \a_j^2 }
\janlk.
\label{01.15}
\qqq
There are no zero self-linking numbers in the step-by-step procedure
applied to
$J^{(\epsilon)}_{\!\aan}(\cL;K)$
for general values of $\epsilon_j$. This means that the analytic
expression for the ``regularized'' invariant
$Z^{({\rm tr})}_\epsilon(M;K)$
obtained through the step-by-step procedure in the assumption that no
self-linking numbers are zero, coincides with the expression coming
from \ex{1.13}. At the same time, since the matrix $l_{ij}$ is
non-degenerate, \ex{1.13} indicates that
$Z^{({\rm tr})}_\epsilon(M;K)$
is a continuous function at $\epsilon_j=0$.
Therefore we can perform a
step-by-step calculation while keeping $\epsilon_j$ as variables and
take a limit $\epsilon_j\rightarrow 0$ only at the very end when all
the integrals have been calculated. The limit exists and is equal to
the \rhs of \ex{1.13}. The possible singularities of intermediate
invariants $J^{({\rm tr})}$ at $\epsilon_j=0$ are ultimately canceled.

%++++++++++++++++++++++++++++++++++++++++++
\nsection{Invariance under the Kirby Meridian Move}
\label{s5}
%++++++++++++++++++++++++++++++++++++++++++
The step-by-step procedure allows for a seemingly simple proof of the
invariance of $\snm$ under the Kirby meridian move. We recall that we
have to prove that performing a surgery on a 1-framed meridian of a
local cable is equivalent to adding an extra
negative
twist to that cable. We
will follow the idea of the Reshetikhin-Turaev proof\cx{RT} of the
invariance of the \wrt invariant as given by the surgery
formula\rx{1.10}. They use the ``fusion'' properties of a local
cable.

Let $J_{\aan}$ be the colored Jones polynomial of a framed link with
an $N$-strand local cable, $\a_j$ being the colors of the cable
strands. We denote the Jones polynomial of the same link with an
extra negative
twist added to the local cable as $\jtw$. $\jmer$ denotes the
Jones polynomial of the link if a 0-framed meridian of the local
cable is added as an extra component carrying the color $\b$. In
accordance with the surgery formula\rx{1.10} Reshetikhin and
Turaev had to show that
\qq
\jtw = e^{-i\pi{3\over 4}{K-2\over K} }
\sum_{1\leq\b\leq K-1} e^{\iptk (\b^2-1) } \jmer \sqtk \sinb.
\label{01.16}
\qqq
If we apply the step-by-step formula\rx{01.3} to the surgery on the
meridian, then we see that the equation
\qq
\jtw = e^{-i\pi{3\over 4}{K-2\over K} }
\sqrt{2K} \int\limits_{0\atop [b=0]}^\infty db\,
e^{\ipkt b^2 -\iptk} \jmer \sin(\pi b)
\label{01.17}
\qqq
would demonstrate the invariance of $\ztrmk$.

Reshetikhin and Turaev used the fact that the polynomial $J_{\aan}$
of a link with a local cable can be decomposed into a sum over the
color $\a$ of a strand into which the local cable is fused:
\qq
J_{\aan} = \sman \jf.
\label{01.18}
\qqq
This decomposition satisfies two properties:
\qq
\jtw & = & \sman \jf e^{-\iptk (\a^2 - 1) },
\label{01.19}
\\
\jmer & = & \sman \jf {\sinab \over \spka}.
\label{01.20}
\qqq
A substitution of these expressions into \eex{01.16} and \rx{01.17}
reveals that it is enough to prove the relations
\qq
e^{-\iptk(\a^2 - 1) } & = & -i^{1\over 2} e^{i\pi\over K}
\sqtk {1\over \spka} \sum_{1\leq \b\leq K-1} e^{\iptk\b^2}
\sinab\sinb,
\label{01.21}\\
e^{-\iptk(\a^2 - 1) } & = & -i^{1\over 2} e^{i\pi\over K}
{\sqrt{2K}\over \spka} \int\limits_{0\atop [b=0]}^\infty db\,
e^{\ipkt b^2} \sin(\pi\a b) \sin(\pi b)
\label{01.22}
\qqq
for \eex{01.16} and \rx{01.17} respectively. Note that the point
$b=0$ is indeed the stationary phase point of the integral in
\ex{01.22} because $\a\ll K$ in view of the bound
$1\leq \a \leq \a_1 + \cdots \a_N$ coming from the summation range in
\ex{01.18}.

Equation\rx{01.22} can be checked by a straightforward calculation
while \ex{01.21} might require an application of the Poisson
resummation. It follows from comparing \eex{01.21} and \rx{01.22}
that the identity
\qq
\sum_{-K\leq \b\leq K-1} e^{\iptk\b^2}
= \int_{-\infty}^{+\infty} d\b\,e^{\iptk\b^2}
\label{01.23}
\qqq
is the reason why the substitution of a stationary phase integral
instead of a sum in the Reshetikhin-Turaev surgery formula\rx{1.10}
does not destroy the invariance of its \lhs under the meridian Kirby
move.

%++++++++++++++++++++++++++++++++++++++
\nsection{Other Results}
\label{s6}
%++++++++++++++++++++++++++++++++++++++

In this final section we want to list the results of \cx{Ro4} and
\cx{Ro5} which are validated by the Propositions~\ref{p1.1},
\ref{p01.1} and by the invariance of $\ztrmk$ under the Kirby move.

A link $\cL$ is called algebraically split (\asl) if all of its
non-diagonal linking numbers are zero:
$l_{ij}=0$ for $i\neq j$. The next proposition is an easy corollary
of the Proposition~\ref{1.1}:
%+++++++++++++++++++++++++++++++
\begin{proposition}
\label{p6.1}
Let $\cL$ be an $N$-component \asl in $S^3$. Then its colored Jones
polynomial has the following $1/K$ expansion (\cf \ex{01.1}):
\qq
\exp \left( -\iptk \sjn l_{jj} \a_j^2 \right) \janlk =
\pjna \sum_{n\geq 0\atop 0\leq m\leq {3\over 4} n}
\Dmn\aaq K^{-n}.
\label{6.1}
\qqq
The coefficients of the polynomials $D_{3n,4n}\aaq$ are expressed in
terms of triple Milnor linking numbers $\lmu_{ijk}$ with the help of
closed diagrams with trivalent vertices.
\end{proposition}
%+++++++++++++++++++++++++++++++++
Although the bound $m\leq {3\over 4}n$ in \ex{6.1} is weaker than the
Melvin-Morton bound $m\leq {n\over 2}$ of \ex{01.1}, still it is
better than the trivial bound $m\leq n$. This allows us to perform
the integrals over $da_j$ of the step-by-step procedure in one scoop:
%+++++++++++++++++++++++++
\begin{proposition}
\label{p6.2}
If $\cL$ is an \asl in $S^3$, then
\qq
\lefteqn{
Z^{({\rm tr})}(\chi_\cL(S^3);K) =
e^{i\ffr} (2K)^{N\over 2} \int\limits_{0\atop [a_j=0]}^\infty
da_1\cdots da_N\, e^{\ipkt \sjn l_{jj} a_j^2}
}\hspace*{2.5in}
\label{6.2}\\
&&\times
\left( e^{-\ipkt\sjn l_{jj} a_j^2} J_{Ka_1,\ldots,Ka_N}(\cL;K)
\right) \pjn \sin(\pi a_j).
\nonumber
\qqq
\end{proposition}
%+++++++++++++++++++++++++++++++++++++++++++

T.~Ohtsuki defined in \cx{Oh1} the finite type invariants of integer
homology spheres. Here is the trivial extension of his definition to
rational homology spheres. Let $\#\cL$ denote a number of components
of a link $\cL\in S^3$. If $\lambda(M)$ is an invariant of \rhs then
we denote by $\tilde{\lambda}(\cL)$ an associated invariant of \asl
$\cL$ with non-zero self-linking numbers defined by the formula
\qq
\tilde{\lambda}(\cL) = \sum_{\cL\p\in \cL}
(-1)^{\#\cL\p} \lambda(\chi_{\cL\p}(S^3) ).
\label{6.3}
\qqq
\begin{definition}
An invariant $\lambda$ of \Rhs is of finite type of at most order $n$
if $\tilde{\lambda}(\cL)=0$ for any \asl $\cL$ with $\#\cL = n+1$. An
invariant $\lambda$ is of order $n$ if it is of at most order $n$ and
not of at most order $n-1$.
\end{definition}
%++++++++++++++++++++
\begin{proposition}
\label{p6.3}
The invariants $\snm$
defined by \eex{1.12},\rx{1.13} are of finite type $3n$.
\end{proposition}
The proof presented in \cx{Ro4} is based on the fact that \ex{6.2}
expresses $S_{n_0}(M)$ in terms of the polynomials $\Dmn\aaq$ with
$n-m\leq n_0$. In view of the bound $m\leq {3\over 4}n$ in the sum of
the \rhs of \ex{6.1} this inequality implies that only the
polynomials $\Dmn$ with $m\leq 3n_0$ participate in the expression
for $S_{n_0}(M)$. Each monomial of the polynomial $\Dmn\aaq$ can
depend on at most $m$ different colors $\a_j$. It turns out that if a
monomial does not depend on a particular color $\a_j$, then its
contribution to $Z^{({\rm tr})}(\chi_{\cL\p};K)$ does not depend on
whether the sublink $\cL\p\in\cL$ contains the link component
$\cL_j$. As a result, the contribution of $\Dmn$ is canceled in the
alternating sum\rx{6.3} if $\#\cL>m$ because for each monomial of
$\Dmn$ there will be a component of $\cL$ whose color is not
represented in that monomial. This proves that $\snm$ is of at most
order $3n$. The diagrammatic rules for the calculation of
$\tilde{S}_n(\cL)$ for $\#\cL = 3n$ developed in \cx{Ro4} demonstrate
that there is a $3n$-component link $\cL$ for which
$\tilde{S}_n(\cL)\neq 0$. This proves that $\snm$ is of exactly order
$3n$.

In their most recent preprint \cx{GO},
S.~Garoufalidis and T.~Ohtsuki gave a
modified definition of the finite type invariants of \Rhs. In
addition to the requirement that
$\tilde{\lambda}(\cL)=0$ for any \asl $\cL$ with $\#\cL = n+1$
they also demanded that $\tilde{\lambda}(\cL)$ for $\#\cL = n$
should have a specific dependence on the self-linking numbers
$l_{jj}$, namely, it should be inversely proportional to the product
$\prod_{j=1}^n l_{jj}$.
It is easy to see that $S_n(M)$ satisfy
this extra requirement. Indeed, the most color diverse monomials of
$D_{3n,4n}$ which contribute to
$\tilde{S_n}(\cL)$ for $\#\cL = 3n$ contain all the colors in the
minimal square power. As a result, the integral in \ex{6.2} produces
exactly the factors of $l_{jj}^{-1}$.

An alternative to look at the structure of $\zmk$ as
$K\rightarrow\infty$ is to study this invariant for prime values of
$K$. This program was developed by R.~Kirby and P.~Melvin\cx{KM},
S.~Garoufalidis\cx{Ga}, L.~Jeffrey\cx{Je}, H.~Murakami\cx{Mu1},
\cx{Mu2}, T.~Ohtsuki\cx{Oh1}, \cx{Oh2}, R.~Lawrence\cx{Lw}, X-S.~Lin
and Z.~Wang\cx{LW}, G.~Masbaum and J.~Roberts\cx{MR},
and
others.
%+++++++++++++++++++++++++++++++++++++
\begin{proposition}[Murakami]
\label{p6.4}
If $K$ is prime, then the modified \wrt invariant of a \Rhs $M$
\qq
Z\p(M;K) =
\left\{
\begin{array}{cl}
{\zmk\over Z(M;3)}, &\mbox{if $K=-1 \mod{4}$}\\
{\zmk\over\overline{Z(M;3)}}, &\mbox{if $K=1 \mod{4}$},
\end{array}
\right.
\label{6.4}
\qqq
defined by Kirby and Melvin\cx{KM} ($Z\p(M;K) = \zmk$ if $M$ is an
integer homology sphere) belongs to the cyclotomic ring
$\ZZ[e^{2\pi i\over K}]$.
\end{proposition}
For $p,q\in\ZZ$, $q\not\equiv 0 \mod{K}$ define
$(p/q)^{\vee}=pq^{\ast}$, here $qq^{\ast}\equiv 1 \mod{K}$.
%+++++++++++++++++++++
\begin{proposition}[Ohtsuki]
\label{p6.5}
For a \Rhs $M$,
let $K$ be a prime number such that
$$\ordH\not\equiv 0\mod{K}.$$
Then
there exists an infinite sequence of rational invariants of $M$:
$\lambda_0(M),\lambda_1(M),\ldots$ such that if we present $Z\p(M;K)$
as a polynomial in $(e^{2\pi i\over K} - 1)$:
\qq
Z\p(M;K) =
\sum_{n=0}^{K-2} a_n (e^{2\pi i\over K} - 1)^n,
\label{6.5}
\qqq
then
\qq
a_n\equiv (\lambda_n)^\vee \mod{K}\qquad\mbox{for}\qquad
0\leq n\leq {K-3\over 2}.
\label{6.6}
\qqq
\end{proposition}
%+++++++++++++++++++++++++++++++++++++++++++++++
In\cx{Ro5} we used the bound $m\leq {3\over 4}n$ in the
expansion\rx{6.1} of the Jones polynomial of an \asl in order to
provide a conceptually simple proof of the Propositions~\ref{p6.4}
and \ref{p6.5}. We also established a relation between the
``perturbative'' invariants $\snm$ and Ohtsuki's invariants
$\lambda_n(M)$:
%++++++++++++++++++++++++++++++++++++
\begin{proposition}
\label{p6.6}
The generating function $\ztrmk$ is proportional to Ohtsuki's
polynomial
\qq
\sum_{n=0}^\infty \lambda_n(M) (e^{2\pi i\over K} - 1)^n =
{\ztrmk\over\zsk}.
\label{6.7}
\qqq
\end{proposition}
%+++++++++++++++++++++++++++++++++++++
Equation\rx{6.7} should be understood in the following way:
for any $n_0>0$
\qq
\sum_{n=0}^{n_0} \lambda_n(M) (e^{2\pi i\over K} - 1)^n
= \frac{\left({\pi\over K} \right)} {\spk}
\exp\left[ \sum_{n=1}^{n_0} \snm \ipk^n \right]
+ \cO(K^{-n_0 - 1}).
\label{6.8}
\qqq
The proof of\cx{Ro5}
is based on an apparent similarity between a gaussian
integral and a gaussian sum for prime $K$:
\qq
\int_{-\infty}^{+\infty}
e^{{2\pi i\over K} \left({p\over q}\a^2 + 2n\a\right)} d\a =
e^{i{\pi\over 4}\sign{{p\over q}} }
\left( {K\over 2} \left|{q\over p}\right| \right)^{1\over 2}
e^{-{2\pi i\over K} {q\over p} n^2},
\label{6.9}
\\
\sum_{\a=0}^{K-1} e^{{2\pi i\over K} (pq^{\ast}\a^2 + 2n\a)}
= e^{i{\pi\over 4} (1-\kappa)} K^{1\over 2}
\left({pq^{\ast}\over K}\right)
e^{-{2\pi i\over K} p^{\ast}qn^2},
\label{6.10}
\qqq
here $p,q,n\in\ZZ$,
\qq
\kappa =
\left\{
\begin{array}{cl}
1 &\mbox{if $K\equiv -1 \mod{4}$}\\
-1&\mbox{if $K\equiv 1 \mod{4}$},
\end{array}
\right.
\label{6.11}
\qqq
and $\left({pq^{\ast}\over K}\right)$ in \ex{6.10} is Legendre's
symbol.

\section*{Acknowledgements}
I am very thankful to the organizers and participants of the School
on Quantum Invariants of 3-manifolds for a wonderful opportunity to
discuss these subjects. I especially appreciate the numerous comments
made by J.~Andersen, S.~Axelrod, D.~Bar-Natan, S.~Garoufalidis,
G.~Masbaum, D.~Thurston and V.~Turaev.

This work was supported by the National Science Foundation
under Grants No. PHY-92 09978 and DMS 9304580. Part of this work was
done during my visit to Argonne National Laboratory. I want to thank
Professor C.~Zachos for his hospitality and advice.

\nappendixe

We are going to comment on the proof of Propositions\rw{p1.1}
and\rw{p01.1} in\cx{Ro2}. We will sketch how this proof can be made
rigorous by the use of Kontsevich's integral. A detailed explanation
will be provided in\cx{Ro9}.

Let $\cL$ be an $N$-dimensional link $\cL$ in $S^3$. Let its
components $\cL_j$ be fixed by smooth maps $x_j(t)$, $0\leq t\leq 1$:
\qq
x_j:\;\;[0,1] \rightarrow S^3, \qquad x_j(0) = x_j(1).
\label{A.1}
\qqq
We also introduce the following notation: $J(x)$ is a 2-form on $S^3$
taking values in the Lie (co-)algebra $su(2)$ (since $su(2)$ has a
standard Killing form, we will not distinguish between the algebra and
its conjugate). In order to simplify combinatorics
in our calculations we will use
variational derivatives in the simplest way possible: if $G(x)$ is a
1-form taking values in $su(2)$, then
\qq
{\delta \over \delta J(x) } \int_{S^3} d^3 y\, G(y)J(y) = G(x)
\label{A.2}
\qqq
(we dropped the wedge: $G(y)J(y) = G(y) \wedge J(y)$, and assumed an
implicit $su(2)$ contraction between $G(y)$ and $J(y)$). According to
\ex{A.2}, $\delta\over \delta J(x) $ behaves as a 1-form on $S^3$,
so we can pull it back with the map\rx{A.1} and integrate it
along the link components.

The proof of the Proposition\rw{p1.1} in\cx{Ro2} was based
essentially on the following assumption about the expansion\rx{1.1}
(unfortunately, we did not emphasize this fact in\cx{Ro2}):
%Now we can formulate the basic property of the expansion\rx{1.1}:
%%%%%%%%%%%%
\begin{assumption}
\label{pA.1}
There exist `multi-local' 1-forms (\ie 1-forms in all their
arguments) $\gmny$, $m\geq 2,\;n\geq0$ on $S^3$ taking values in
$(su(2))^{\otimes m}$ such that the expansion\rx{1.1} comes from the
formula
\qq
\janlk & = &
%\snjz \!\pjn \;\itanj \Tr_{\a_j} \ddjpanj
\pjn \left[ \sum_{n_j\geq 0} \;\itanj \Tr_{\a_j} \ddjpanj \right]
\label{A.4}
\\
&&
%\qquad
\times
\ekiy{\smtnz}.
\nonumber
\qqq
This equation should be understood in the following way: for any
$n_0 > 0$
\qq
\lefteqn{
\pjna \sum_{0\leq n \leq n_0 \atop 0\leq m\leq n} \Dmn\aaq K^{-n}
}
\label{A.5}\\
& = &
%\sum_{ {n_j\geq0 \atop (1\leq j\leq N)} \atop {\sjN n_j \leq 2n_0} }
%\pjn \; \itanj \Tr_{\a_j} \ddjpanj
\pjn \left[ \sum_{n_j \geq 0 \atop \sjN n_j \leq 2n_0}
\;\itanj \Tr_{\a_j} \ddjpanj \right]
\nonumber\\
&&\qquad\times
\ekiy{\smtnzs}
\nonumber\\
&&
\qquad\qquad
+ \cO(K^{-n_0 - 1} )
\nonumber
\qqq
\end{assumption}
%%%%%%%%%%%%%%%%%%%%%%%%%%%%%%
It is clear from \ex{A.5} that we treat the sum $\smtnz$ in the
exponential of the \rhs of \ex{A.4} as a formal power series in
$K^{-1}$. We do not need to know whether it is convergent, because
only a finite number of terms in this sum contribute to any
particular term in the \rhs of \ex{1.1}.

%The Proposition\rw{pA.1}
Equation\rx{A.4}
is the standard assumption of the quantum
field theory applied to the Chern-Simons action in the framework of
Witten's description\cx{Wi1} of the colored Jones
polynomial. The set of the forms $\gmny$ is not unique, it depends
on the choice of gauge fixing required for the calculation of the
Chern-Simons path integral.

M.~Kontsevich proved the Assumption\rw{pA.1} (modulo some minor
corrections, see
\eg\cx{BN2} and references therein) without any reference to path
integrals. He obtained a particularly simple set of forms $G_{m,n}$:
%
%The choice of the forms $\gmny$ in \ex{A.4} is not unique.
%M.~Kontsevich derived \ex{A.4}  with
%
\qq
\gmny = 0 \qquad \mbox{if either}\;\;m\geq 3 \;\;\mbox{or}\;\;
n\geq 1,
\label{A0.5}
\qqq
so that only $G_{2,0} \neq 0$. The approach of D.~Bar-Natan\cx{BN1}
and R.~Bott, C.~Taubes\cx{BT} produces different forms $G_{m,n}$.

%In\cx{Ro2} we derived the Proposition\rw{p1.1} essentially from
%\ex{A.4} (unfortunately,
%we did not emphasize this fact there).
%this fact was not emphasized there).

The
particular form of the forms $G_{m,n}$ played no role
in the proof of the Proposition\rw{p1.1} in\cx{Ro2},
because that
proof involved only manipulations with path ordered integrals
$$\itanj$$
and Lie algebra traces $\Tr_{\a_j}$.
Indeed, the proof of\cx{Ro2}
was based on the following general fact which can
be derived from the Campbell-Hausdorf formula
(\cf Proposition~2.2
of\cx{Ro2}, there we used the notation $A_\mu(x)$ for the 1-form
${\delta\over \delta J(x) }$):
%%%%%%%%%%%%%%%%%
\begin{proposition}
\label{pA.2}
There exist $SU(2)$-equivariant multilinear forms $\cten{n}(\vdot)$,
\linebreak
$\vdot \in su(2)$ with values in $su(2)$:
\qq
\cten{n}:\;\;(su(2))^{\otimes n} \rightarrow su(2)
\label{A.6}
\qqq
(in particular,
$\cten{1}(v) = v$, $\cten{2}(v_1,v_2) = {1\over 2} [v_1,v_2]$ ), such
that
\qq
& \sum_{n_j\geq 0}\; \itanj\, \ddja{n_j} = e^{v_j},
\label{A.7}\\
& v_j = \sum_{n\geq 1}\; \itan\; \cten{n} \ddjpan,
\;\;\;\; v_j \in su(2).
\label{A.8}
\qqq
\end{proposition}
Another ingredient in the proof of\cx{Ro2} was Kirillov's trace
formula
\qq
\Tr_{\a_j} e^{v_j} =
{|\vv_j|\over \sin|\vv_j|} \int_{|\vaa_j| = \a_j}
{d^2\vaa_j\over 4\pi \a_j}\, e^{i \vv_j\cdot \vaa_j}.
\label{A.9}
\qqq
We use vector notations for the elements of $su(2)$ in the \rhs of
this formula in order to emphasize that $su(2)$ is isomorphic to
$\IR^3$, $\cdot$ is the Killing form and $|\cdot|$ is the
corresponding norm.
The combination of eqs.\rx{A.4},\rx{A.7},\rx{A.8} and\rx{A.9} leads
to the Reshetikhin's formula\rx{1.2} through the application of
rigorous combinatorial rules (which are associated with Feynman rules
in quantum field theory) to the calculation of the action of
derivatives $\delta\over \delta J$ contained in
$e^{\vv_j\cdot \vaa_j}$, on the exponential in the \rhs of \ex{A.4}.
We will provide the details in\cx{Ro9}.

The proof of the Proposition\rw{p01.1} in\cx{Ro2} was based
on two assumptions about the properties of the forms $G_{m,0}$,
$m\geq 2$.
%%%%%%
\begin{assumption}
\label{Aa1}
For a bi-local $(1,1)$-form
\qq
\Omega(y_1,y_2) = {2\over i\pi} G_{2,0}(y_1,y_2),
\label{A.10}
\qqq
there exists a bi-local $(0,2)$-form $\tilde{\Omega}(y_1,y_2)$
such that
\qq
d_{y_2} \Omega(y_1,y_2) = \delta^{(3)} (y_2 - y_1) +
d_{y_1} \tilde{\Omega}(y_1,y_2).
\label{A.13}
\qqq
\end{assumption}
%%%%%%%%%%%%
\begin{assumption}
\label{Aa2}
The forms $G_{m,0}$, $m\geq 3$ come from tree Feynman diagrams with
propagator\rx{A.10} and triple vertices of the Chern-Simons action
(see, \eg\cx{BN1} and references therein for
details).
\end{assumption}
%%%%%%%%%%

In\cx{Ro2} we derived the set of forms $G_{m,n}$ with the help of the
quantum field theory perturbation theory described \eg in\cx{BN1}.
Therefore Assumption\rw{Aa2} was satisfied automatically. Also, it
was explained in\cx{BN1} that
in the notations $S^3= \IR^3 \bigcup \{\infty\}$, $y_{1,2}\in \IR^3$
the form\rx{A.10} was equal to
\qq
\Omega(y_1,y_2) = \ytof\, dy_1^\mu\,dy_2^\nu,
\label{A.11}
\qqq
and the Assumption\rw{Aa1} was satisfied with the choice of the form
\qq
\tilde{\Omega}(y_1,y_2) = \ytof\, dy_2^\mu \wedge dy_2^\nu.
\label{A.12}
\qqq

The two Assumptions\rw{Aa1} and\rw{Aa2} are also satisfied for the set
of forms $G_{m,0}$ coming from Kontsevich integral derivation of
\ex{A.4}.
We use the
notations
$S^3 = \IC^1 \times \IR^1\;\bigcup\;\{\infty\}$,
$(z,t)\in \IC^1\times \IR^1$, so Kontsevich's form $\Omega$ is
\qq
\Omega\zzbt = \dtz\, (dt_1 dz_2 - dt_2 dz_1).
\label{A.14}
\qqq
It satisfies Assumption\rw{Aa1} if we choose
\qq
\tilde{\Omega}\zzbt = -\dtz\, dt_2\wedge dz_2.
\label{A.15}
\qqq
The Feynman diagrams with cubic vertices built upon the
propagator\rx{A.14} are all zero, because the form\rx{A.14} does not
contain $d\bar{z}$ so that its triple wedge products are zero.
Therefore Assumption\rw{Aa2} is satisfied in the trivial way due
to \ex{A0.5}.
Hence the proof of Proposition\rw{p01.1} remains valid
{\em verbatim} if we use the set of forms $G_{m,n}$ coming
from Kontsevich's integral instead of the ones coming from the
\cx{BN1}-style perturbation theory.

%%%%%%%%%%%%%%%%%%%%%%%%%%%%%%%%%%
%%%%%%%%%%%%%%%%%%%%%%%%%%%%%%%%%%

\end{document}